\documentclass[12pt,a4paper,leqno]{amsart}
\usepackage[english]{babel}
\usepackage[T1]{fontenc}
\usepackage{amssymb}
\usepackage{amsfonts}
\usepackage{amsmath}
\usepackage{mathptm}
\swapnumbers
\title{\large{Removability results for subharmonic  functions, for separately subharmonic functions, for harmonic functions, for separately harmonic functions  and for holomorphic functions, a survey}}
\author{\large{Juhani Riihentaus}}
\begin{document}
\null{}
\vskip7cm
\maketitle
\vspace{1ex}
\begin{center}{\large{University of Eastern Finland, Department of Physics and Mathematics,}}\end{center}
\begin{center}{\large{P.O. Box 111, 
FI-80101 Joensuu,  
Finland}}\end{center}

\vspace{0.3ex}

\begin{center}{\large{and}}\end{center}

\vspace{0.3ex}

\begin{center}{\large{University of Oulu,
Department of Mathematical Sciences,}}\end{center}
\begin{center}{\large{P.O. Box 3000,
FI-90014 Oulun yliopisto,
Finland}}\end{center}

\vspace{3ex}

\centerline{{\large{juhani.riihentaus@gmail.com}} \,\, {\large{and}} \,\, {\large{{juhani.riihentaus{@}uef.fi}}}}

\vspace{3ex}
\centerline{\large{ORCID 0000-0002-3631-762X}}
\thispagestyle{empty}
\newpage

\date{\noindent 21 May, 2018}
\large{
\markright{\large{Removable singularities for ${\mathcal{C}}^2$  subharmonic functions}}
\maketitle
\vskip0.3cm
\setcounter{page}{2}

\vspace*{0.3cm}
\begin{abstract}
Blanchet  has shown that a ${\mathcal{C}}^2$ subharmonic  function can be extended through a ${\mathcal{C}}^1$ hypersurface provided the function is continuous throughout and satisfies certain ${\mathcal{C}}^1$-type continuity conditions on the exceptional hypersurface. Later we improved Blanchet's result, at least in a certain sense, by measuring the exceptional set with the aid of  Hausdorff measure.  Our result was flexible, and with the aid of it, we gave certain extension results for harmonic and for holomorphic functions, related to Besicovitch's and Shiffman's well-known extension results, at least in some sense. Now we return to this subject and give a unified presentation of the until now obtained results. First, we refine our  subharmonic function extension result slightly still more, improving also our previous proof. Though our result might be considered a little bit technical and even complicated, it is, nevertheless, flexible. As  examples of its flexibility, we recall here our recent concise extension result for subharmonic functions. Moreover, we slightly refine our previous corollaries for harmonic and for holomorphic functions. In addition, and as  another application, we recall our recent extension result for separately subharmonic functions. And as a new result, we give a related extension result for separately  harmonic functions. As another topic, we recall a slightly related extension result for holomorphic functions.
\end{abstract}
\footnotetext[1]{\noindent 2000 \emph{Mathematics Subject Classification.} 31B05, 31B25, 32A10, 32D20.}
\footnotetext[2]{\noindent \emph{Key words and phrases.} Subharmonic function, separately subharmonic function, harmonic function, separately harmonic function, holomorphic function, Hausdorff measure, exceptional sets}
\setcounter{equation}{0}
\section{Introduction}
\subsection{An outline}
We will consider extension problems for subharmonic, separately subharmonic,  harmonic, separately harmonic and for holomorphic functions. Our results will be based on an extension result for subharmonic functions, see Theorem~1 below in Section~2,  and also \cite{Ri19}, Theorem~1, p.~61. As a matter of fact, our Theorem~1 is a refined version of our previous  result, see \cite{Ri15}, Theorem, p.~568, and \cite{Ri16$_2$}, Theorem~1, p.~154. With the aid of this previous version of Theorem~1, we have already given certain extension results for  harmonic and for holomorphic functions, see \cite{Ri16$_1$}, Theorem~3, p.~51, and \cite{Ri16$_2$}, Theorem~2, p.~155, Theorem~3, p.~156. Using now this improved version of Theorem~1,  we slightly improve these results in  Theorem~3 and in Theorem~5 below in Sections~4 and 6. 

With the aid of our Theorem~1, we have also given an extension result for separately subharmonic functions, see \cite{Ri19}, Theorem~2, p.~62, and Corollary~2,  p.~64. Now we slightly improve the proof of this result in Theorem~2 below in  Section~3. Moreover, and as a new result, we give an extension result for separately harmonic functions, see Theorem~4 and Corollary~4 below in Section 5.  As an additional new result, we give in Theorem~6 below in Section~6, an extension result for holomorphic functions, which is slightly related to our holomorphic function extension results Theorem~5 and Corollary~5.

Observe that the starting point for our results and especially for our Theorem~1, is a result of Blanchet \cite{Bl95}, Theorem~3.1, p. 312. As a matter of fact, in his result Blanchet has shown that hypersurfaces of class ${\mathcal{C}}^1$ are removable singularities for subharmonic functions, provided the considered subharmonic  functions satisfy certain extra assumptions. 
\subsection{Notation}
Our notation is more or less standard, see \cite{{Ri04},{Ri13},{Ri15},{Ri16$_1$},{Ri16$_2$},{Ri19}}. However, for the convenience of the reader we recall here the following.
We use the common convention $0\cdot \pm \infty =0$. For each $n\geq 1$ we identify ${\mathbb{C}}^n$ with ${\mathbb{R}}^{2n}$.  In integrals we write  $dx$ or $dm_n$ for the Lebesgue measure in ${\mathbb{R}}^{n}$, $n\in {\mathbb{N}}$. 
Let $0\leq \alpha \leq n$ and $A\subset {\mathbb{R}}^n$, 
$n\geq 1$. Then we write ${\mathcal{H}}^{\alpha}(A)$ for the $\alpha$-dimensional Hausdorff (outer) measure of $A$. Recall that $\mathcal{H}^0(A)$ is the number of points of $A$.  If $x=(x_1,x_2,\dots ,x_n)\in {\mathbb R}^n$, $n\geq 2$, and $j\in {\mathbb{N}}$, $1\leq j\leq n$, then we write $x=(x_j,X_j)$, where 
$X_j=(x_1,x_2,\dots ,x_{j-1},x_{j+1},\dots ,x_n)$. Moreover, if $A\subset {\mathbb{R}}^n$, $1\leq j\leq n$, and 
$x_j^0\in {\mathbb{R}}$, 
$X_j^0\in {\mathbb{R}}^{n-1}$, we write
\[ A(x_j^0)=\{\, X_j\in {\mathbb{R}}^{n-1}\, :\, \, x=(x_j^0,X_j)\in A\, \},\, \, \, 
 A(X_j^0)=\{\, x_j\in {\mathbb{R}}\, :\, \, x=(x_j,X_j^0)\in A\, \}.\]
If  $\Omega \subset {\mathbb{R}}^n$ and $p>0$, then
  ${\mathcal L}^p_{\mathrm{loc}}(\Omega )$, $p>0$, is 
the space of
functions $u$ in $\Omega$ for which  $\vert u\vert^p$ is locally integrable on $\Omega$.      

For the definition and properties of harmonic and subharmonic functions,  see e.g. \cite{{ArGa01},{HaKe76},{Hel69},{Her71},{Le69},{Ra37}}, see also \cite{{Ri07},{Ri09}}. For the definition and properties of holomorphic functions see e.g. \cite{{Ch89},{JaPf00},{JaPf11}}.   
\section{Extension results for subharmonic functions}
\subsection{A result of Federer}
The following important result of Federer on geometric measure theory will be used repeatedly.

\vskip0.3cm

\noindent \textbf{Lemma.} (\cite{Fe69},  Theorem 2.10.25, p. 188, and \cite{Sh68}, Corollary 4, Lemma~2, p.~114) {\emph{Suppose that $E\subset {\mathbb{R}}^n$, $n\geq 2$. Let $\alpha \geq 0$ and let $\pi_k :{\mathbb{R}}^n \to {\mathbb{R}}^k$ denote the projection onto the first $k$ coordinates.
\begin{itemize}
\item[{1.}] If 
 ${\mathcal{H}}^{k+\alpha}(E)=0$, then ${\mathcal{H}}^{\alpha}(E\cap \pi_k ^{-1}\{x\})=0$ for ${\mathcal{H}}^{k}$-almost all $x\in {\mathbb{R}}^{k}$.
\item[{2.}] If ${\mathcal{H}}^{k+\alpha}(E)<+\infty$, then  ${\mathcal{H}}^{\alpha}(E\cap \pi_k ^{-1}\{x\})<+\infty$ for ${\mathcal{H}}^{k}$-almost all $x\in {\mathbb{R}}^{k}$. 
\end{itemize}}}
\subsection{A result of Blanchet} Blanchet has given the following  result: 

\vskip0.3cm

\noindent \textbf{Blanchet's theorem.} (\cite{Bl95}, Theorems 3.1, 3.2 and 3.3, pp. 312-313) {\emph{Let
 $\Omega$ be a domain in ${\mathbb{R}}^n$, $n\geq 2$, and let $S$ be a hypersurface of class ${\mathcal{C}}^1$ which 
divides $\Omega$ into  two subdomains $\Omega_1$ and $\Omega_2$. Let $u\in {\mathcal{C}}^0(\Omega )\cap 
{\mathcal{C}}^2(\Omega_1\cup \Omega_2)$ be subharmonic (respectively convex (or respectively plurisubharmonic 
provided $\Omega$ is then a domain in ${\mathbb{C}}^n$, $n\geq 1$)) in $\Omega_1$ and $\Omega_2$. If $u_i=u\vert \Omega_i\in {\mathcal{C}}^1
(\Omega_i\cup S)$,
 $i=1,2$, and
\begin{equation} \frac{\partial u_i}{\partial \overline{n}^k}\geq \frac{\partial u_k}{\partial \overline{n}^k}\end{equation}
on $S$ with $i,k=1,2$, then $u$ is subharmonic (respectively convex (or respectively plurisubharmonic)) in $\Omega$.}}

Above $\overline{n}^k=(\overline{n}_1^k,\overline{n}_2^k, \dots ,\overline{n}_n^k)$ is the unit normal 
exterior to $\Omega_k$, and $u_k\in {\mathcal{C}}^1(\Omega_k\cup S)$, $k=1,2$, means that there exist $n$ functions $v^j_k$,
 $j=1,2, \dots ,n$, continuous on $\Omega_k\cup S$, such that 
\[ v_k^j(x)=\frac{\partial u_k}{\partial x_j}(x)\]
for all  $x\in \Omega_k$,  $k=1,2$ and  $j=1,2, \dots ,n$.

\vskip0.3cm

The following example shows that  one cannot drop the above condition (1) in Blanchet's theorem, see also \cite{{Ri15},{Ri16$_2$}}.

\vskip0.3cm  

\noindent \textbf{Example~1.} The function  $u:{\mathbb{R}}^2 \rightarrow {\mathbb{R}}$,
\begin{displaymath}u(z)=u(x+iy)=u(x,y):=\left\{ \begin{array}{ll}1+x, &\textrm{when $x<0$},\\
1-x, &\textrm{when $x\geq 0$,}
\end{array}\right.
\end{displaymath}
is continuous in ${\mathbb{R}}^2$ and subharmonic, even harmonic in ${\mathbb{R}}^2\setminus (\{0\}\times {\mathbb{R}})$. It is easy to see that $u$ does not satisfy the condition (1) on $S=\{0\}\times {\mathbb{R}}$ and that $u$ is not subharmonic in ${\mathbb{R}}^2$.

\vskip0.3cm 

\noindent \textbf{Remark~1.} For related results, previous and later, see Khabibullin's results \cite{Kh92}, Lemma~2.2, p.~201,  Fundamental Theorem~2.1, pp.~200-201, and \cite{Kh99}, Lemma~4.1, p.~503, Theorem~2.1, p.~498, Theorems~3.1 and 3.2, pp.~500-501. In this connection, see also \cite{HaKe76}, 1.4.3, pp.~21-22.
\subsection{An improvement to the result of Blanchet}
Already in \cite{Ri04}, Theorem~4, pp.~181-182, we have given partial improvements to the cited subharmonic  removability results of  Blanchet. For other previous  improvements, see \cite{Ri15}, Theorem, p.~568, and \cite{Ri16$_2$}, Theorem~1, p.~154. Instead of hypersurfaces of class ${\mathcal{C}}^1$, we have there considered arbitrary sets of finite \mbox{$(n-1)$}-dimensional Hausdorff measure as exceptional sets.  However, now we recall, and use,  our most recent result,  \cite{Ri19}, Theorem~1, p.~61, see Theorem~1 below, where  Blanchet's condition $(1)$ is replaced by   the condition (iv), which is now, at least seemingly, less stringent than in  our previous results  \cite{Ri15}, Theorem, (iv), p.~568, and \cite{Ri16$_2$}, Theorem~1, (iv), p.~154, say. 

\vskip0.3cm

\noindent \textbf{Theorem~1.} (\cite{Ri19}, Theorem~1, p.~61)
{\emph{Suppose that  $\Omega$ is a domain in ${\mathbb{R}}^n$, $n\geq 2$. 
Let $E\subset \Omega$ be closed in $\Omega$ 
and let ${\mathcal{H}}^{n-1}(E)<+\infty$. Let 
$u:\Omega \setminus E \to {\mathbb{R}}$ be subharmonic and such that the following conditions are satisfied:
\begin{itemize}
\item[(i)] $u\in {\mathcal{L}}_{{\mathrm{loc}}}^1(\Omega )$.
\item[(ii)] $u\in {\mathcal{C}}^2(\Omega \setminus E)$.
\item[(iii)] For each $j$, $1\leq j\leq n$, $\frac{\partial ^2u}{\partial x_j^2}\in {\mathcal{L}}_{{\mathrm{loc}}}^1(\Omega )$.
\item[(iv)] For each $j$, $1\leq j\leq n$,  and for ${\mathcal{H}}^{n-1}$-almost all $X_j\in {\mathbb{R}}^{n-1}$   such that $E(X_j)$ is finite, the following condition holds:\newline
For each $x_j^0\in E(X_j)$ there exist sequences $x_{j,l}^{0,1}, x_{j,l}^{0,2}\in (\Omega \setminus E)(X_j)$, $l=1,2,\dots$, such that $x_{j,l}^{0,1}\nearrow x_{j}^{0}$, $x_{j,l}^{0,2}\searrow x_{j}^{0}$ as $l\to +\infty$, and 
\begin{center}
\begin{itemize}
\item[(iv(a))]  
$\lim_{l\to +\infty}u(x_{j,l}^{0,1},X_j)=\lim_{l\to +\infty}u(x_{j,l}^{0,2},X_j)\in {\mathbb{R}},$
\item[(iv(b))] $-\infty <\liminf_{l\to +\infty}\frac{\partial u}{\partial x_j}(x_{j,l}^{0,1},X_j)\leq \limsup_{l\to +\infty} 
\frac{\partial u}{\partial x_j}(x_{j,l}^{0,2},X_j)<+\infty .$\
\end{itemize}
\end{center}
 \end{itemize}
Then $u$ has a subharmonic extension to $\Omega$.}}}

\vskip0.3cm

\noindent{\textit{Proof.}} Observe first that using suitable subsequences one can replace the assumption (iv) by the following, (only) seemingly stronger condition:
\vskip0.3cm
\begin{itemize}
\item[(iv$^*$)]
\noindent \emph{For each $x_j^0\in E(X_j)$ there exist sequences $x_{j,l}^{0,1}, x_{j,l}^{0,2}\in (\Omega \setminus E)(X_j)$, $l=1,2,\dots$, such that $x_{j,l}^{0,1}\nearrow x_{j}^{0}$, $x_{j,l}^{0,2}\searrow x_{j}^{0}$ as $l\to +\infty$, and} 
\begin{center}
\begin{itemize}
\item[\,\,\,\,(iv$^*$(a))]  
$\lim_{l\to +\infty}u(x_{j,l}^{0,1},X_j)=\lim_{l\to +\infty}u(x_{j,l}^{0,2},X_j)\in {\mathbb{R}},$
\item[\,\,\,\,(iv$^*$(b))] $-\infty <\lim_{l\to +\infty}\frac{\partial u}{\partial x_j}(x_{j,l}^{0,1},X_j)\leq \lim_{l\to +\infty} 
\frac{\partial u}{\partial x_j}(x_{j,l}^{0,2},X_j)<+\infty .$\
\end{itemize}
\end{center}
\end{itemize}

\vskip0.3cm

It is sufficient to show that 
\begin{equation*}\int u(x)\, \varDelta \varphi (x)\, dx\geq 0\end{equation*}
for all nonnegative testfunctions $\varphi \in {\mathcal{D}}(\Omega )$, see e.g. \cite{Her71}, Corollary~1, p.~13.

Take $\varphi \in {\mathcal{D}}(\Omega )$, $\varphi \geq 0$, arbitrarily. Let $K=\mathrm{spt}\,\varphi$.  Choose a domain $\Omega_1$ such that $K\subset \Omega_1\subset \overline{\Omega}_1\subset \Omega$ 
and $\overline{\Omega}_1$ is compact. Since $u\in {\mathcal{C}}^2(\Omega \setminus E)$
and $u$ is subharmonic in $\Omega \setminus E$, $\varDelta u(x)\geq 0$ for all $x\in \Omega \setminus E$. Thus the claim  follows if we show that 
\begin{equation*}\int u(x)\, \varDelta \varphi (x)\, dx\geq \int  \varDelta u(x)\, \varphi (x)\, dx.\end{equation*}

For this purpose fix $j$, $1\leq j\leq n$, arbitrarily for a while. By Fubini's theorem,  see e.g. \cite{EvGa93}, Theorem~1,  pp.~22-23,
\begin{equation*}\int u(x)\, \frac{\partial ^2\varphi}{\partial x_j^2}(x)\, dx
= \int \bigl[\int u(x_j,X_j)\frac{\partial ^2 \varphi}{\partial x_j^2}(x_j,X_j)\, dx_j\bigr]\, dX_j.\end{equation*}
Using the above Lemma, assumptions (i), (ii) and (iii), and Fubini's theorem, we see that for ${\mathcal{H}}^{n-1}$-almost all
 $X_j\in {\mathbb{R}}^{n-1}$,
\begin{equation} \left\{ \begin{aligned}
&u( \cdot \, \,  ,X_j)\in {\mathcal{L}}_{\mathrm{loc}}^1(\Omega (X_j)),\\
&\frac{\partial^2 u}{\partial x_j^2}(\, \,  \cdot \, \,  ,X_j)\in {\mathcal{L}}_{\mathrm{loc}}^1(\Omega (X_j)),\\
&E(X_j){\textrm{ \, is finite, thus there exists \,}} M=M(X_j)\in {\mathbb{N}}_0 {\mathrm{\, \, such\, \, that \,\,}}\\
&E(X_j)=\{\, x_j^1,x_j^2, \ldots ,x_j^M\, \} {\mathrm{\,\, where\,\,}} x_j^k<x_j^{k+1},\,\,k=1,2, \ldots ,M-1.\\
\end{aligned}\right.
\end{equation}
Let $X_j\in {\mathbb{R}}^{n-1}$ be arbitrary as above in (2). We may suppose that $\Omega (X_j)$ is a finite interval. Choose for each $k=1,2, \dots ,M$ numbers $a_k, b_k\in (\Omega \setminus E)(X_j)$ such that 
$a_k<x_j^k<b_k, \,\, k=1,2,\dots ,M$, $a_{k+1}=b_k$, $k=1,2,\dots, M-1$, and that $a_1,b_M\in (\Omega\setminus \overline{\Omega}_1)(X_j)$.

With the aid of (iv$^*$) we find for each $x_j^k\in E(X_j)$  sequences $x_{j,l
}^{k,1}, x_{j,l}^{k,2}\in (\Omega \setminus E)(X_j)$, $l=1,2, \dots$, $x_{j,l}^{k,1}\nearrow x_{j}^{k}$, $x_{j,l}^{k,2}\searrow x_{j}^{k}$ as $l\to +\infty$, such that 
\begin{equation*}\lim_{l\to +\infty} u(x_{j,l}^{k,1},X_j)=\lim_{l\to +\infty}u(x_{j,l}^{k,2},X_j)\in {\mathbb{R}},
\end{equation*}
and
\begin{equation*}-\infty <\lim_{l\to +\infty}\frac{\partial u}{\partial x_j}(x_{j,l}^{k,1},X_j)\leq \lim_{l\to +\infty}
\frac{\partial u}{\partial x_j}(x_{j,l}^{k,2},X_j)<+\infty .
\end{equation*}
 Take $k$, $1\leq k\leq M$, arbitrarily and  consider the  interval $(a_k,b_k)$, where $a_k<x_j^k<b_k$. To simplify the notation, write $a:=a_k$, $b:=b_k$ and $x_j^0:=x_j^k$.
 Then
 \[a<x_{j,l}^{0,1}\nearrow x_{j}^{0}, \,\,\, b>x_{j,l}^{0,2}\searrow x_{j}^{0}\,\,\,  {\mathrm{ as }} \,\,\, l\rightarrow +\infty .\]  
 
Using partial integration we get: 
\begin{equation*}
\begin{split}
\int_a^b &u(x_j,X_j) \, \frac{\partial ^2\varphi}{\partial x_j^2}(x_j,X_j)\, dx_j=
\int_a^{x_j^0} u(x_j,X_j)\, \frac{\partial ^2\varphi}{\partial x_j^2}(x_j,X_j)\, dx_j+\\&+
\int_{x_j^0}^b u(x_j,X_j)\, \frac{\partial ^2\varphi}{\partial x_j^2}(x_j,X_j)\, dx_j\\
&=\lim_{l\rightarrow +\infty}\int_a^{x_{j,l}^{0,1}} u(x_j,X_j)\, \frac{\partial ^2\varphi}{\partial x_j^2}(x_j,X_j)\, dx_j+
\lim_{l\rightarrow +\infty}\int_{x_{j,l}^{0,2}}^b u(x_j,X_j)\, \frac{\partial ^2\varphi}{\partial x_j^2}(x_j,X_j)\, dx_j\\
&=\lim_{l\rightarrow +\infty}\left[{\Big{\vert}}_a^{x_{j,l}^{0,1}}u(x_j,X_j)\, \frac{\partial \varphi}{\partial x_j}(x_j,X_j)- \int_a^{x_{j,l}^{0,1}} \frac{\partial u}{\partial x_j}(x_j,X_j)\, \frac{\partial \varphi}{\partial x_j}(x_j,X_j)\, dx_j \right]+\\
&+\lim_{l\rightarrow +\infty}\left[{\Big{\vert}}_{x_{j,l}^{0,2}}^b\,u(x_j,X_j)\frac{\partial \varphi}{\partial x_j}(x_j,X_j)-\int_{x_{j,l}^{0,2}}^b \frac{\partial u}{\partial x_j}(x_j,X_j)\, \frac{\partial \varphi}{\partial x_j}(x_j,X_j)\, dx_j\right] \\
&=\left[u(b,X_j)\frac{\partial \varphi}{\partial x_j}(b,X_j)-u(a,X_j)\frac{\partial \varphi}{\partial x_j}(a,X_j)\right]+\\
&+\lim_{l\rightarrow +\infty}\left[u({x_{j,l}^{0,1}},X_j)\, \frac{\partial \varphi}{\partial x_j}(x_{j,l}^{0,1},X_j)- \int_a^{x_{j,l}^{0,1}}\frac{\partial u}{\partial x_j} (x_j,X_j)\, \frac{\partial \varphi}{\partial x_j}(x_j,X_j)\, dx_j \right]+\\
&-\lim_{l\rightarrow +\infty}\left[u({x_{j,l}^{0,2}},X_j)\, \frac{\partial \varphi}{\partial x_j}(x_{j,l}^{0,2},X_j)+ \int_{x_{j,l}^{0,2}}^b\frac{\partial u}{\partial x_j} (x_j,X_j)\, \frac{\partial \varphi}{\partial x_j}(x_j,X_j)\, dx_j \right]\\
&=\left[u(b,X_j)\frac{\partial \varphi}{\partial x_j}(b,X_j)-u(a,X_j)\frac{\partial \varphi}{\partial x_j}(a,X_j)\right]+\\
&-\lim_{l\rightarrow +\infty} \int_a^{x_{j,l}^{0,1}}\frac{\partial u}{\partial x_j} (x_j,X_j)\, \frac{\partial \varphi}{\partial x_j}(x_j,X_j)\, dx_j +\\
&-\lim_{l\rightarrow +\infty}\int_{x_{j,l}^{0,2}}^b\frac{\partial u}{\partial x_j} (x_j,X_j)\, \frac{\partial \varphi}{\partial x_j}(x_j,X_j)\, dx_j \\
&=\left[u(b,X_j)\frac{\partial \varphi}{\partial x_j}(b,X_j)-u(a,X_j)\frac{\partial \varphi}{\partial x_j}(a,X_j)\right]+\\
&-\lim_{l\rightarrow +\infty}\left[{\Big{\vert}}_a^{x_{j,l}^{0,1}} \frac{\partial u}{\partial x_j}(x_j,X_j)\varphi (x_j,X_j)- \int_a^{x_{j,l}^{0,1}}\frac{\partial^2 u}{\partial x_j^2} (x_j,X_j)\, \varphi (x_j,X_j)\, dx_j \right]+\\
&-\lim_{l\rightarrow +\infty}\left[{\Big{\vert}}_{x_{j,l}^{0,2}}^b \frac{\partial u}{\partial x_j}(x_j,X_j)\varphi (x_j,X_j)- \int_{x_{j,l}^{0,2}}^b\frac{\partial^2 u}{\partial x_j^2} (x_j,X_j)\, \varphi (x_j,X_j)\, dx_j \right]\\
\end{split}
\end{equation*}
\begin{equation*}
\begin{split}
&=\left[u(b,X_j)\frac{\partial \varphi}{\partial x_j}(b,X_j)-u(a,X_j)\frac{\partial \varphi}{\partial x_j}(a,X_j)\right]+\\
&+\left[\frac{\partial u}{\partial x_j}(a,X_j)\varphi (a,X_j)-\frac{\partial u}{\partial x_j}(b,X_j)\varphi (b,X_j)\right]+\\
&-\lim_{l\rightarrow +\infty}\left[\frac{\partial u}{\partial x_j}(x_{j,l}^{0,1},X_j)\varphi (x_{j,l}^{0,1},X_j)\right]+\\
&+\lim_{l\rightarrow +\infty}\left[\frac{\partial u}{\partial x_j}(x_{j,l}^{0,2},X_j)\varphi (x_{j,l}^{0,2},X_j)\right]+
\int_a^b \frac{\partial^2 u}{\partial x_j^2}(x_j,X_j)\, \varphi (x_j,X_j)\, dx_j \\
&=\left[u(b,X_j)\frac{\partial \varphi}{\partial x_j}(b,X_j)-u(a,X_j)\frac{\partial \varphi}{\partial x_j}(a,X_j)\right]+\\
&+\left[\frac{\partial u}{\partial x_j}(a,X_j)\varphi (a,X_j)-\frac{\partial u}{\partial x_j}(b,X_j)\varphi (b,X_j)\right]+\\
&+\left[\lim_{l\rightarrow +\infty}\frac{\partial u}{\partial x_j}(x_{j,l}^{0,2},X_j)-\lim_{l\rightarrow +\infty}\frac{\partial u}{\partial x_j}(x_{j,l}^{0,1},X_j)\right]\varphi (x_j^0,X_j)+\\
&+\int_a^b \frac{\partial^2 u}{\partial x_j^2}(x_j,X_j)\, \varphi (x_j,X_j)\, dx_j \\
&\geq \left[u(b,X_j)\frac{\partial \varphi}{\partial x_j}(b,X_j)-u(a,X_j)\frac{\partial \varphi}{\partial x_j}(a,X_j)\right]+\\
&+\left[\frac{\partial u}{\partial x_j}(a,X_j)\varphi (a,X_j)-\frac{\partial u}{\partial x_j}(b,X_j)\varphi (b,X_j)\right]+\int_a^b \frac{\partial^2 u}{\partial x_j^2}(x_j,X_j)\, \varphi (x_j,X_j)\, dx_j . 
\end{split}
\end{equation*}
Above we have used just standard properties of limits and our assumption (iv$^*$(b)).  

To return to our original notation, we have thus obtained for each $k=1,2, \dots ,M$,
\begin{equation*}\begin{split}
&\int_{a_k}^{b_k} u(x_j,X_j) \, \frac{\partial ^2\varphi}{\partial x_j^2}(x_j,X_j)\, dx_j\geq
 \left[u(b_k,X_j)\frac{\partial \varphi}{\partial x_j}(b_k,X_j)-u(a_k,X_j)\frac{\partial \varphi}{\partial x_j}(a_k,X_j)\right]+\\
&+\left[\frac{\partial u}{\partial x_j}(a_k,X_j)\varphi (a_k,X_j)-\frac{\partial u}{\partial x_j}(b_k,X_j)\varphi (b_k,X_j)\right]+\int_{a_k}^{b_k} \frac{\partial^2 u}{\partial x_j^2}(x_j,X_j)\, \varphi (x_j,X_j)\, dx_j. \\
\end{split}\end{equation*}
Then just sum over $k=1,2, \dots ,M$:
\begin{equation*}\begin{split}
&\int  u(x_j,X_j) \, \frac{\partial ^2\varphi}{\partial x_j^2}(x_j,X_j)\, dx_j=
\int_{a_1}^{b_M} u(x_j,X_j) \, \frac{\partial ^2\varphi}{\partial x_j^2}(x_j,X_j)\, dx_j=\\
&=\sum_{k=1}^M\int_{a_k}^{b_k} u(x_j,X_j) \, \frac{\partial ^2\varphi}{\partial x_j^2}(x_j,X_j)\, dx_j\geq \\
\end{split}\end{equation*}
\begin{equation*}\begin{split}
&\geq \sum_{k=1}^M \left[u(b_k,X_j)\frac{\partial \varphi}{\partial x_j}(b_k,X_j)-u(a_k,X_j)\frac{\partial \varphi}{\partial x_j}(a_k,X_j)\right]+\\
&+\sum_{k=1}^M\left[\frac{\partial u}{\partial x_j}(a_k,X_j)\varphi (a_k,X_j)-\frac{\partial u}{\partial x_j}(b_k,X_j)\varphi (b_k,X_j)\right]+\\
&+\sum_{k=1}^M\int_{a_k}^{b_k} \frac{\partial^2 u}{\partial x_j^2}(x_j,X_j)\, \varphi (x_j,X_j)\, dx_j = \\
&=\int_{a_1}^{b_M} \frac{\partial^2 u}{\partial x_j^2}(x_j,X_j)\, \varphi (x_j,X_j)\, dx_j=
\int \frac{\partial^2 u}{\partial x_j^2}(x_j,X_j)\, \varphi (x_j,X_j)\, dx_j .  
\end{split}\end{equation*}
Above we have used the choice of the numbers $a_k,b_k$, $k=1,2, \dots ,M$, and the fact that $a_1,b_M\in (\Omega\setminus \overline{\Omega}_1)(X_j)$. 

Integrate then with respect to $X_j$ and use again Fubini's theorem:
\begin{equation*}\begin{split}
\int &u(x)\frac{\partial^2 \varphi}{\partial x_j^2}(x)\,dx=\int \left[\int  u(x_j,X_j) \, \frac{\partial ^2\varphi}{\partial x_j^2}(x_j,X_j)\, dx_j\right]\,dX_j\geq \\
&\geq \int \left[ \int  \frac{\partial ^2
u}{\partial x_j^2}(x_j,X_j)\varphi (x_j,X_j)\, dx_j\right]dX_j=\int \frac{\partial ^2u}{\partial x_j^2}(x)\varphi (x)\,dx.  
\end{split}\end{equation*}
Summing over $j=1,2, \dots ,n$ gives the desired inequality
\begin{equation*}
\begin{split}
\int & u(x)\varDelta \varphi(x)\,dx=\int  u(x) \,\sum_{j=1}^n \frac{\partial^2 \varphi}{\partial x_j^2}(x)\,dx\geq \\
&\geq \int   \,\sum_{j=1}^n \frac{\partial^2 u}{\partial x_j^2}(x)\varphi (x)\,dx=\int \varDelta u(x)\varphi (x)\,dx \geq 0,\\
\end{split}
\end{equation*}
concluding the proof. \hfill \qed
\vskip0.3cm
\noindent \textbf{Example~2.} The function  $u:{\mathbb{R}}^2 \rightarrow {\mathbb{R}}$, given already in Example~1, is continuous in ${\mathbb{R}}^2$ and subharmonic, even harmonic in ${\mathbb{R}}^2\setminus (\{0\}\times {\mathbb{R}})$, but not subharmonic in ${\mathbb{R}}^2$. Observe that $u$ satisfies the above conditions (i), (ii), (iii) and (iv(a)) in ${\mathbb{R}}^2\setminus (\{0\}\times {\mathbb{R}})$. However, $u\vert {\mathbb{R}}^2\setminus (\{0\}\times {\mathbb{R}})$ does not satisfy the condition (iv(b)). Thus one cannot drop the condition (iv(b)) in Theorem~1.  
\subsection{ } The following corollary reflects the strength of Theorem~1. 

\vskip0.3cm

\noindent \textbf{Corollary~1.} (\cite{Ri19}, Corollary~1, p.~61) {\emph{Suppose that  $\Omega$ is a domain in ${\mathbb{R}}^n$, $n\geq 2$. Let $E\subset \Omega$ be closed in $\Omega$ 
and let ${\mathcal{H}}^{n-1}(E)=0$. Let 
$u:\Omega \setminus E \to {\mathbb{R}}$ be subharmonic and such that the following conditions are satisfied:
\begin{itemize}
\item[(i)] $u\in {\mathcal{L}}_{\mathrm{loc}}^1(\Omega )$,
\item[(ii)] $u\in {\mathcal{C}}^2(\Omega \setminus E)$,
\item[(iii)] for each $j$, $1\leq j\leq n$, 
{\Large{$\frac{\partial ^2u}{\partial x_j^2}$}}$\in {\mathcal{L}}_{{\mathrm{loc}}}^1
(\Omega )$.
\end{itemize}
Then $u$ has a subharmonic extension to $\Omega$.}}

\vskip0.3cm

\noindent{\textit{Proof.}} Follows directly from Theorem~1 and from the above Lemma of  Federer.
\hfill \qed
\section{Extension results for separately subharmonic functions}
\subsection{} Next we will give an extension result for separately subharmonic functions. Our proof will be based on Theorem~1 and on the following nice result, Proposition~1 below. Observe  that the there used \emph{hypoharmonic functions} are in our terminology just \emph{subharmonic functions}.  

\vskip0.3cm

\noindent \textbf{Proposition~1.} (\cite{Her71}, Proposition~1, p.~33)
{\emph{Suppose that  $\Omega$ is a domain in ${\mathbb{R}}^{p+q}$, $p, q\geq 2$. 
Let $w:\Omega \to [-\infty ,+\infty )$ be nearly subharmonic. Let $w^*:\Omega \to [-\infty ,+\infty )$ be the regularized function of $w$, which is then subharmonic. Then the following properties are equivalent.
\begin {itemize}
\item[(1)] The distribution $\Delta_x w=\Delta_x w^*$= (sum of the square second order derivatives of $w$ or $w^*$ with respect to the $p$ coordinates of $x$) is positive. 
\item[(2)] For all $y\in {\mathbb{R}}^q$ the function 
\begin{center} 
$\Omega (y)\ni x\mapsto w^*(x,y)\in [-\infty ,+\infty )$
\end{center}
is hypoharmonic.
\item[(3)] For almost every $y\in {\mathbb{R}}^q$ the function 
\begin{center} 
$\Omega(y)\ni x\mapsto w^*(x,y)\in [-\infty ,+\infty )$
\end{center}
is subharmonic.
\item[(4)] For almost every $y\in {\mathbb{R}}^q$ the function 
\begin{center} 
$\Omega(y)\ni x\mapsto w(x,y)\in [-\infty ,+\infty )$
\end{center}
is nearly subharmonic.
\end{itemize}}}

\vskip0.3cm

\subsection{An extension result for separately subharmonic functions} Then our result. Observe that our proof below is a slightly improved version of our  original proof in \cite{Ri19}, pp.~62-64.

\vskip0.3cm

\noindent \textbf{Theorem~2.} (\cite{Ri19}, Theorem~2, p.~62)
{\emph{Suppose that  $\Omega$ is a domain in ${\mathbb{R}}^{p+q}$, $p, q\geq 2$. 
Let $E\subset \Omega$ be closed in $\Omega$ 
and let ${\mathcal{H}}^{p+q-1}(E)<+\infty$. Let 
$w:\Omega \setminus E \to {\mathbb{R}}$ be separately subharmonic, that is, 
\begin{center} 
for all $y\in {\mathbb{R}}^q$ the function $(\Omega\setminus E)(y)\ni x\mapsto w(x,y)\in {\mathbb{R}}$ is subharmonic,
\end{center}
and
\begin{center} 
for all $x\in {\mathbb{R}}^p$ the function $(\Omega\setminus E)(x)\ni y\mapsto w(x,y)\in {\mathbb{R}}$ is subharmonic,
\end{center}
and such that the following conditions are satisfied:
\begin{itemize}
\item[(i)] $w\in {\mathcal{L}}_{{\mathrm{loc}}}^1(\Omega )$.
\item[(ii)] $w\in {\mathcal{C}}^2(\Omega \setminus E)$.
\item[(iii)] For each $j$, $1\leq j\leq p$, $\frac{\partial ^2w}{\partial x_j^2}\in {\mathcal{L}}_{{\mathrm{loc}}}^1(\Omega )$, and for each $k$, $1\leq k\leq q$, $\frac{\partial ^2w}{\partial y_k^2}\in {\mathcal{L}}_{{\mathrm{loc}}}^1(\Omega )$.
\item[(iv)] For each $j$, $1\leq j\leq p$,  and for ${\mathcal{H}}^{p-1+q}$-almost all $(X_j,y)\in {\mathbb{R}}^{p-1+q}$   such that $E(X_j,y)$ is finite, the following condition holds:\newline
For each $x_j^0\in E(X_j,y)$ there exist sequences $x_{j,l}^{0,1}, x_{j,l}^{0,2}\in (\Omega \setminus E)(X_j,y)$, $l=1,2,\dots$, such that $x_{j,l}^{0,1}\nearrow x_{j}^{0}$, $x_{j,l}^{0,2}\searrow x_{j}^{0}$ as $l\to +\infty$, and 
\begin{center}
\begin{itemize}
\item[(iv(a))]  
$\lim_{l\to +\infty}w(x_{j,l}^{0,1},X_j,y)=\lim_{l\to +\infty}w(x_{j,l}^{0,2},X_j,y)\in {\mathbb{R}},$
\item[(iv(b))] $-\infty <\liminf_{l\to +\infty}\frac{\partial w}{\partial x_j}(x_{j,l}^{0,1},X_j,y)\leq \limsup_{l\to +\infty} 
\frac{\partial w}{\partial x_j}(x_{j,l}^{0,2},X_j,y)<+\infty .$
\end{itemize}
\end{center}
\item[(v)] For each $k$, $1\leq k\leq q$, and for ${\mathcal{H}}^{p+q-1}$-almost all $(x,Y_k)\in {\mathbb{R}}^{p+q-1}$   such that $E(x,Y_k)$ is finite, the following condition holds:\newline
For each $y_k^0\in E(x,Y_k)$ there exist sequences $y_{k,l}^{0,1}, y_{k,l}^{0,2}\in (\Omega \setminus E)(x,Y_k)$, $l=1,2,\dots$, such that $y_{k,l}^{0,1}\nearrow y_{k}^{0}$, $y_{k,l}^{0,2}\searrow y_{k}^{0}$ as $l\to +\infty$, and 
\begin{center}
\begin{itemize}
\item[(v(a))]  
$\lim_{l\to +\infty}w(x,y_{k,l}^{0,1},Y_k)=\lim_{l\to +\infty}w(x,y_{k,l}^{0,2},Y_k)\in {\mathbb{R}},$
\item[(v(b))] $-\infty <\liminf_{l\to +\infty}\frac{\partial w}{\partial y_k}(x,y_{k,l}^{0,1},Y_k)\leq \limsup_{l\to +\infty} 
\frac{\partial w}{\partial y_k}(x,y_{k,l}^{0,2},Y_k)<+\infty .$
\end{itemize}
\end{center} 
 \end{itemize}
Then $w$ has a separately subharmonic extension to $\Omega$.}}

\vskip0.3cm

\noindent{\textit{Proof.}} Since $w$ is separately subharmonic in $\Omega\setminus E$ and $w\in {\mathcal{L}}_{\textrm{loc}}^1(\Omega)$, it follows from \cite{Ri07}, Theorem~A, p.~50, or from \cite{Ri13}, Corollary~4.6, p.~412, say, that $w$ is subharmonic in $\Omega\setminus E$. Using then Theorem~1 one sees that $w:\Omega \setminus E \to {\mathbb{R}}$ has a subharmonic extension $w^*:\Omega \to [-\infty ,+\infty )$.  

To show that the subharmonic function $w^*:\Omega \to [-\infty ,+\infty )$ is in fact separately subharmonic, we must show that for all $y\in {\mathbb{R}}^q$ the function 
\begin{center} 
$\Omega(y)\ni x\mapsto w^*(x,y)\in [-\infty,+\infty )$
\end{center}
is subharmonic, and that for all  $x\in {\mathbb{R}}^p$ the function 
\begin{center} 
$\Omega(x)\ni y\mapsto w^*(x,y)\in [-\infty,+\infty )$
\end{center}
is subharmonic. It is clearly sufficient to show that the first claim holds.

For this purpose we define 
${\tilde{w}}:\Omega \rightarrow [-\infty,+\infty )$,
\begin{displaymath}{\tilde{w}}(x,y):=\left\{ \begin{array}{ll}w(x,y), &\textrm{when}\,\, (x,y)\in \Omega\setminus E,\\
-\infty, &\textrm{when}\,\, (x,y)\in E.
\end{array}\right.
\end{displaymath}
To see that $\tilde{w}$ is nearly subharmonic, observe first that 
$\tilde{w}(x,y)=w(x,y)$ for all $(x,y)\in \Omega\setminus E$, thus for almost all $(x,y)\in \Omega$. Hence $\tilde{w}\in {\mathcal{L}}^1_{\textrm{loc}}(\Omega)$. To see that the mean value inequality holds, take $(x_0,y_0)\in \Omega$ and $\overline{B^{p+q}((x_0,y_0),r)}\subset\Omega$ arbitrarily. If $(x_0,y_0)\in \Omega \setminus E$, then 
\begin{align*}
{\tilde{w}}(x_0,y_0)&=w(x_0,y_0)=w^*(x_0,y_0)\leq \\
&\leq \frac{1}{\nu_{p+q}r^{p+q}}\int_{B^{p+q}((x_0,y_0),r)}w^*(x,y)\,dm_{p+q}(x,y)=\\
&\leq \frac{1}{\nu_{p+q}r^{p+q}}\int_{B^{p+q}((x_0,y_0),r)}w(x,y)\,dm_{p+q}(x,y)=\\
&\leq \frac{1}{\nu_{p+q}r^{p+q}}\int_{B^{p+q}((x_0,y_0),r)}{\tilde{w}}(x,y)\,dm_{p+q}(x,y).\\
\end{align*} 
If $(x_0,y_0)\in E$, then, since $\tilde{w}\in {\mathcal{L}}^1_{\mathrm{loc}}(\Omega)$, 
\begin{center}
$-\infty = \tilde{w}(x_0,y_0)< \frac{1}{\nu_{p+q}r^{p+q}}\int_{B^{p+q}((x_0,y_0),r)}{\tilde{w}}(x,y)\,dm_{p+q}(x,y)\in {\mathbb{R}}.$
\end{center}
Thus $\tilde{w}$ is nearly subharmonic.

In order to be able to use Proposition~1, we show 
that for  ${\mathcal{H}}^{q}$-almost all $y\in {\mathbb{R}}^q$, the function 
\begin{center} 
$\Omega(y)\ni x\mapsto \tilde{w}(x,y)\in [-\infty,+\infty )$
\end{center}
is nearly subharmonic. 

For this purpose  fix $j$, $1\leq j\leq p$, arbitrarily for a while. 
 
By our assumption  ${\mathcal{H}}^{p-1+q}(E)<+\infty$. From the above Lemma of Federer it follows that for ${\mathcal{H}}^{p-1+q}$-almost all $(X_j,y)\in {\mathbb{R}}^{p-1+q}$ the set $E(y)(X_j)$ is finite. Write 
\[A_j:=\{(X_j,y)\in {\mathbb{R}}^{p-1+q}\,:\,E(y)(X_j) \,\, {\textrm{is finite}}\}.\]
Thus 
\[{\mathcal{H}}^{p-1+q}(A_j^c)=0\, \Leftrightarrow \, m_{ p-1+q}(A_j^c)=0\,  \Leftrightarrow \, \int_{{\mathbb{R}}^{p-1+q}}\,\chi_{A_j^c}(X_j,y)
dm_{p-1+q}(X_j,y)=0,\] 
where $\chi_{A_j^c}(\cdot ,\cdot )$ is the characteristic function of the set $A_j^c$, the complement taken in ${\mathbb{R}}^{p-1+q}$. 

Next use Fubini's theorem:
\[0= \int_{{\mathbb{R}}^{p-1+q}}\chi_{A_j^c}(X_j,y)dm_{p-1+q}(X_j,y)
=\int_{{\mathbb{R}}^q}[\int_{{\mathbb{R}}^{p-1}}\chi_{A_j^c}(X_j,y)dm_{p-1}(X_j)]dm_q(y).\]
Since 
\[\int_{{\mathbb{R}}^{p-1}}\chi_{A_j^c}(X_j,y)dm_{p-1}(X_j)\geq 0,\]
we see that in fact 
\[\int_{{\mathbb{R}}^{p-1}}\chi_{A_j^c}(X_j,y)dm_{p-1}(X_j)= 0\]
for ${\mathcal{H}}^{q}$-almost all $y\in{\mathbb{R}}^{q}$. 

Write 
\begin{equation*}\begin{split}
B_1^j:&=\{y\in {\mathbb{R}}^q\, :\, \int_{{\mathbb{R}}^{p-1}}\chi_{A_j^c}(X_j,y)dm_{p-1}(X_j)=0\,\},\\
&=\{y\in {\mathbb{R}}^q \, : \, \chi_{A_j^c}(X_j,y)=0 \,\, {\textrm{for}} \,\,{\mathcal{H}}^{p-1}-{\textrm{almost all}}\,\, X_j\in {\mathbb{R}}^{p-1}\},\\
&=\{y\in {\mathbb{R}}^q \, : \, \chi_{A_j}(X_j,y)=1 \,\, {\textrm{for}} \,\,{\mathcal{H}}^{p-1}- {\textrm{almost all}}\,\, X_j\in {\mathbb{R}}^{p-1}\},\\
&=\{y\in {\mathbb{R}}^q \, : \, E(y)(X_j)\,\, {\textrm{is finite for}} \,\,{\mathcal{H}}^{p-1}- {\textrm{almost all}}\,\, X_j\in {\mathbb{R}}^{p-1}\}.\\
\end{split}
\end{equation*}
Write $B_1:=B^1_1\cap B^2_1\cap\dots \cap B^p_1$. Then for all $y\in B_1$ we have $(X_j,y)\in A_j$, that is, $E(y)(X_j)$ is finite for ${\mathcal{H}}^{p-1}$-almost all $X_j\in {\mathbb{R}}^{p-1}$, and this for all $j=1,2,\dots ,p$.

Next write 
\begin{equation*}
\begin{split}
B_2&:=\{y\in {\mathbb{R}}^q\, : \, {\mathcal{H}}^{p-1}(E(y))<+\infty \,\},\\ 
B_3&:=\{y\in {\mathbb{R}}^q \, : \, w(\cdot ,y)\in {\mathcal{L}}^1_{{\textrm{loc}}}(\Omega(y))\}, \\
B_4&:=\{y\in {\mathbb{R}}^q \, : \, w(\cdot ,y)\in {\mathcal{C}}^2((\Omega \setminus E)(y))\},\\
B_5^j&:=\{y\in {\mathbb{R}}^q \, : \, \frac{\partial^2}{\partial x_j^2}w(\cdot ,y)\in {\mathcal{L}}^1_{{\textrm{loc}}}(\Omega(y))\},\\
B_5&:=B_5^1\cap B_5^2\cap \dots \cap B_5^p,\\
B:&=B_1\cap B_2 \cap B_3 \cap B_4\cap B_5.\\
\end{split}
\end{equation*}
 
Then for all  $y\in B$ the function 
\[(\Omega\setminus E)(y)\ni x\mapsto w(x,y)\in {\mathbb{R}}\]
satisfies the assumptions of Theorem~1. Therefore these functions  have nearly subharmonic extensions
\begin{center} 
$\Omega(y)\ni x\mapsto w^{**}(x,y)\in [-\infty ,+\infty )$.
\end{center}

To complete the  proof, we show that for all $y\in B$ the function 
\[(\Omega)(y)\ni x\mapsto \tilde{w}(x,y)\in [-\infty ,+\infty )\]
is nearly subharmonic. Observe first that $\tilde{w}(\cdot ,y)\in {\mathcal{L}}_{\mathrm{loc}}^1(\Omega(y))$, since $y\in B$. To show  that $\tilde{w}(\cdot ,y)$, $y\in B$, satisfies the mean value inequality, take $x_0\in (\Omega \setminus E)(y)$ arbitrarily. Since ${\mathcal{H}}^{p-1}(E(y))<+ \infty$, we have
\begin{equation*}
\begin{split}
\tilde{w}(x_0,y)=& w(x_0,y)=w^{**}(x_0,y)\leq \frac{1}{\nu_pr^p}\int_{B^p(x_0,r)}\, w^{**}(x,y)\,dm_p(x)=\\
\leq & \frac{1}{\nu_p r^p}\int_{B^p(x_0,r)}\, w(x,y)\,dm_p(x)= \frac{1}{\nu_p r^p}\int_{B^p(x_0,r)}\tilde{w}(x,y)dm_p(x).\\
\end{split}
\end{equation*}
In the case $x_0\in E(y)$ we have $\tilde{w}(x_0,y)=-\infty$, and thus the mean value inequality is automatically satisfied.

Therefore our claim follows from Proposition~1. \hfill \qed

\vskip0.3cm

\noindent \textbf{Example~3.} (\cite{Ri19}, Example~1, p.~64) The function  $u:{\mathbb{R}}^4 \rightarrow {\mathbb{R}}$,
\begin{displaymath}u(z_1,z_2)=u(x_1+iy_1,x_2+iy_2)=u(x_1,y_1,x_2,y_2):=\left\{ \begin{array}{ll}1+x_1, &\textrm{when $x_1<0$},\\
1-x_1, &\textrm{when $x_1\geq 0$,}
\end{array}\right.
\end{displaymath}
is continuous in ${\mathbb{R}}^4$ and separately subharmonic, even separately harmonic in ${\mathbb{R}}^4\setminus (\{0\}\times {\mathbb{R}}^3)$, but not separately subharmonic in ${\mathbb{R}}^4$. As a matter of fact, it is easy to see that $u$ is not even subharmonic in ${\mathbb{R}}^4$.  Observe that $u$ satisfies the above conditions (i), (ii), (iii), (iv(a)) and (v(a)) in ${\mathbb{R}}^4\setminus (\{0\}\times {\mathbb{R}}^3)$. However, $u\vert {\mathbb{R}}^4\setminus (\{0\}\times {\mathbb{R}}^3)$ does not satisfy the conditions (iv(b)) and (v(b)). Thus these conditions cannot be dropped in Theorem~2.

\vskip0.3cm

\noindent \textbf{Corollary~2.} (\cite{Ri19}, Corollary~2, p.~64) {\emph{Suppose that  $\Omega$ is a domain in ${\mathbb{R}}^{p+q}$, $p, q\geq 2$. Let $E\subset \Omega$ be closed in $\Omega$ 
and let ${\mathcal{H}}^{p+q-1}(E)=0$. Let 
$w:\Omega \setminus E \to {\mathbb{R}}$ be separately subharmonic, that is, 
\begin{center} 
$for \,\,all\,\, y\in {\mathbb{R}}^q\,\, the \,\, function\,\,(\Omega\setminus E)(y)\ni x\mapsto w(x,y)\in {\mathbb{R}} \,\,is\,\, subharmonic$,
\end{center}
and
\begin{center}
$for \,\,all\,\, x\in {\mathbb{R}}^p\,\, the \,\, function\,\,(\Omega\setminus E)(x)\ni y\mapsto w(x,y)\in {\mathbb{R}} \,\,is\,\, subharmonic$, 
\end{center}
and such that the following conditions are satisfied:
\begin{itemize}
\item[(i)] $w\in {\mathcal{L}}_{\mathrm{loc}}^1(\Omega )$,
\item[(ii)] $w\in {\mathcal{C}}^2(\Omega \setminus E)$,
\item[(iii)] for each $j$, $1\leq j\leq p$, 
{\Large{$\frac{\partial ^2w}{\partial x_j^2}$}}$\in {\mathcal{L}}_{{\mathrm{loc}}}^1
(\Omega )$ and for each $k$, $1\leq k\leq q$, 
{\Large{$\frac{\partial ^2w}{\partial y_k^2}$}}$\in {\mathcal{L}}_{{\mathrm{loc}}}^1
(\Omega )$.
\end{itemize}
Then $w$ has a separately  subharmonic extension to $\Omega$.}}

\vskip0.3cm

\noindent{\textit{Proof.}} Follows directly from Theorem~2 and from the above Lemma of  Federer.
\hfill \qed

\section{Extension results for harmonic functions}
\subsection{}
For removability results for harmonic functions see,  among others,  \cite{{HaPo70},{HaPo75},{Po84},{Ul91}} and the references therein, say.

Now, using  our Theorem~1, we give the following extension result for harmonic functions, compare with our previous version \cite{Ri16$_2$}, Theorem~2, p.~155:

\vskip0.3cm

\noindent \textbf{Theorem~3.}
{\emph{Suppose that  $\Omega$ is a domain in ${\mathbb{R}}^n$, $n\geq 2$. 
Let $E\subset \Omega$ be closed in $\Omega$ 
and let ${\mathcal{H}}^{n-1}(E)<+\infty$. Let 
$u:\Omega \setminus E\to {\mathbb{R}}$ be harmonic and such that the following conditions are satisfied:
\begin{itemize}
\item[(i)] $u\in {\mathcal{L}}_{{\mathrm{loc}}}^1(\Omega )$.
\item[(ii)] For each $j$, $1\leq j\leq n$, $\frac{\partial ^2u}{\partial x_j^2}\in {\mathcal{L}}_{{\mathrm{loc}}}^1(\Omega )$.
\item[(iii)] For each $j$, $1\leq j\leq n$,  and for ${\mathcal{H}}^{n-1}$-almost all $X_j\in {\mathbb{R}}^{n-1}$   such that $E(X_j)$ is finite, the following condition holds:\newline
For each $x_j^0\in E(X_j)$ there exist sequences $x_{j,l}^{0,1}, x_{j,l}^{0,2}\in (\Omega \setminus E)(X_j)$, $l=1,2,\dots$, such that $x_{j,l}^{0,1}\nearrow x_{j}^{0}$, $x_{j,l}^{0,2}\searrow x_{j}^{0}$ as $l\to +\infty$, and 
\begin{center}
\begin{itemize}
\item[(iii(a))]  
$\lim_{l\to +\infty}u(x_{j,l}^{0,1},X_j)=\lim_{l\to +\infty}u(x_{j,1}^{0,2},X_j)\in {\mathbb{R}},$
\item[(iii(b))] $-\infty <\liminf_{l\to +\infty}\frac{\partial u}{\partial x_j}(x_{j,l}^{0,1},X_j)= \limsup_{l\to +\infty} 
\frac{\partial u}{\partial x_j}(x_{j,l}^{0,2},X_j)<+\infty .$\
\end{itemize}
\end{center}
\end{itemize}
Then $u$ has a unique harmonic extension to $\Omega$.}}

\vskip0.3cm

\noindent{\textit{Proof.}} Since the assumptions of Theorem~1 do hold for the subharmonic function $u$, $u$ has a subharmonic extension $u^*$ to $\Omega$. On the other hand, the assumptions of Theorem~1 hold also for the subharmonic function $v=-u$. Thus $v=-u$ has a subharmonic extension $v^{*}=(-u)^*$ to $\Omega$. As above in the proof of Theorem~1, we may suppose that the limits
\begin{equation*}
\lim_{l\rightarrow +\infty}\left[\frac{\partial u}{\partial x_j}(x_{j,l}^{0,1},X_j)\right]\,\,\,\, {\mathrm{ and }}\,\, \,\, \lim_{l\rightarrow +\infty}\left[\frac{\partial u}{\partial x_j}(x_{j,l}^{0,2},X_j)\right]
\end{equation*}
indeed exist.

Since $-v^*=u^*$, the extension $u^*$ of $u$ is both subharmonic and superharmonic, thus harmonic and the claim follows.  \hfill \qed 
\vskip0.3cm
\noindent \textbf{Example~4.} The function  $u:{\mathbb{R}}^2 \rightarrow {\mathbb{R}}$ given already in Example~2, shows that one cannot drop the condition (iii(b)) in Theorem~3.  
\subsection{ } Then a concise special case to  Theorem~3:

\vskip0.3cm

\noindent \textbf{Corollary~3.} (\cite{Ri16$_2$}, Corollary~3, p.~155) {\emph{Suppose that  $\Omega$ is a domain in ${\mathbb{R}}^n$, $n\geq 2$. 
Let $E\subset \Omega$ be closed in $\Omega$ 
and let ${\mathcal{H}}^{n-1}(E)=0$. Let 
$u:\Omega \setminus E \to {\mathbb{R}}$ be harmonic and  such that the following conditions are satisfied:
\begin{itemize}
\item[(i)] $u\in {\mathcal{L}}_{{\mathrm{loc}}}^1(\Omega )$,
\item[(ii)] for each $j$, $1\leq j\leq n$, $\frac{\partial ^2u}{\partial x_j^2}\in {\mathcal{L}}_{{\mathrm{loc}}}^1(\Omega )$.
\end{itemize}
Then $u$ has a unique harmonic extension to $\Omega$.}}

\vskip0.3cm

\noindent{\textit{Proof.}} With the aid of the above Lemma one sees easily that the assumptions of  Theorem~3 are satisfied. \hfill \qed 
\section{Extension results for separately harmonic functions}
\subsection{} Next we  give an extension result for separately harmonic functions. Our proof will be based on Theorem~2 and on basic properties of subharmonic functions. 

\vskip0.3cm

\noindent \textbf{Theorem~4.}
{\emph{Suppose that  $\Omega$ is a domain in ${\mathbb{R}}^{p+q}$, $p, q\geq 2$. 
Let $E\subset \Omega$ be closed in $\Omega$ 
and let ${\mathcal{H}}^{p+q-1}(E)<+\infty$. Let 
$w:\Omega \setminus E \to {\mathbb{R}}$ be separately harmonic, that is, 
\begin{center} 
for all $y\in {\mathbb{R}}^q$ the function $(\Omega\setminus E)(y)\ni x\mapsto w(x,y)\in {\mathbb{R}}$ is harmonic,
\end{center}
and
\begin{center} 
for all $x\in {\mathbb{R}}^p$ the function $(\Omega\setminus E)(x)\ni y\mapsto w(x,y)\in {\mathbb{R}}$ is harmonic,
\end{center}
and such that the following conditions are satisfied:
\begin{itemize}
\item[(i)] $w\in {\mathcal{L}}_{{\mathrm{loc}}}^1(\Omega )$.
\item[(ii)] For each $j$, $1\leq j\leq p$, $\frac{\partial ^2w}{\partial x_j^2}\in {\mathcal{L}}_{{\mathrm{loc}}}^1(\Omega )$, and for each $k$, $1\leq k\leq q$, $\frac{\partial ^2w}{\partial y_k^2}\in {\mathcal{L}}_{{\mathrm{loc}}}^1(\Omega )$.
\item[(iii)] For each $j$, $1\leq j\leq p$,  and for ${\mathcal{H}}^{p-1+q}$-almost all $(X_j,y)\in {\mathbb{R}}^{p-1+q}$   such that $E(X_j,y)$ is finite, the following condition holds:\newline
For each $x_j^0\in E(X_j,y)$ there exist sequences $x_{j,l}^{0,1}, x_{j,l}^{0,2}\in (\Omega \setminus E)(X_j,y)$, $l=1,2,\dots$, such that $x_{j,l}^{0,1}\nearrow x_{j}^{0}$, $x_{j,l}^{0,2}\searrow x_{j}^{0}$ as $l\to +\infty$, and 
\begin{center}
\begin{itemize}
\item[(iii(a))]  
$\lim_{l\to +\infty}w(x_{j,l}^{0,1},X_j,y)=\lim_{l\to +\infty}w(x_{j,l}^{0,2},X_j,y)\in {\mathbb{R}},$
\item[(iii(b))] $-\infty <\liminf_{l\to +\infty}\frac{\partial w}{\partial x_j}(x_{j,l}^{0,1},X_j,y)= \limsup_{l\to +\infty} 
\frac{\partial w}{\partial x_j}(x_{j,l}^{0,2},X_j,y)<+\infty .$
\end{itemize}
\end{center}
\item[(iv)] For each $k$, $1\leq k\leq q$, and for ${\mathcal{H}}^{p+q-1}$-almost all $(x,Y_k)\in {\mathbb{R}}^{p+q-1}$   such that $E(x,Y_k)$ is finite, the following condition holds:\newline
For each $y_k^0\in E(x,Y_k)$ there exist sequences $y_{k,l}^{0,1}, y_{k,l}^{0,2}\in (\Omega \setminus E)(x,Y_k)$, $l=1,2,\dots$, such that $y_{k,l}^{0,1}\nearrow y_{k}^{0}$, $y_{k,l}^{0,2}\searrow y_{k}^{0}$ as $l\to +\infty$, and 
\begin{center}
\begin{itemize}
\item[(iv(a))]  
$\lim_{l\to +\infty}w(x,y_{k,l}^{0,1},Y_k)=\lim_{l\to +\infty}w(x,y_{k,l}^{0,2},Y_k)\in {\mathbb{R}},$
\item[(iv(b))] $-\infty <\liminf_{l\to +\infty}\frac{\partial w}{\partial y_k}(x,y_{k,l}^{0,1},Y_k)=\limsup_{l\to +\infty} 
\frac{\partial w}{\partial y_k}(x,y_{k,l}^{0,2},Y_k)<+\infty .$
\end{itemize}
\end{center} 
 \end{itemize}
Then $w$ has a separately harmonic extension to $\Omega$.}}

\vskip0.3cm

\noindent{\textit{Proof.}} By Lelong's theorem, see \cite{Le61} or \cite{Her71}, Theorem, p.~54, separately harmonic functions are  harmonic.  Thus we know that  $w\in {\mathcal{C}}^2(\Omega \setminus E)$. From Theorem~2 it then follows that $w$ has a separately subharmonic extension $w^*: \Omega \to [-\infty ,+\infty )$. Similarly, $u=-w$ has a separately subharmonic extension  $u^*: \Omega \to [-\infty , +\infty)$. It is easy to see that $w^*=-u^*$. Thus $w^*$ is in fact separately harmonic.  \hfill \qed

\vskip0.3cm

\noindent \textbf{Example~5.} The function  $u:{\mathbb{R}}^4 \rightarrow {\mathbb{R}}$ given already in Example~3, shows that one cannot drop the conditions (iii(b)) and (iv(b)) in Theorem~4.

\vskip0.3cm

\noindent \textbf{Corollary~4.} {\emph{Suppose that  $\Omega$ is a domain in ${\mathbb{R}}^{p+q}$, $p, q\geq 2$. Let $E\subset \Omega$ be closed in $\Omega$ 
and let ${\mathcal{H}}^{p+q-1}(E)=0$. Let 
$w:\Omega \setminus E \to {\mathbb{R}}$ be separately harmonic, that is, 
\begin{center} 
$for \,\,all\,\, y\in {\mathbb{R}}^q\,\, the \,\, function\,\,(\Omega\setminus E)(y)\ni x\mapsto w(x,y)\in {\mathbb{R}} \,\,is\,\, harmonic$,
\end{center}
and
\begin{center}
$for \,\,all\,\, x\in {\mathbb{R}}^p\,\, the \,\, function\,\,(\Omega\setminus E)(x)\ni y\mapsto w(x,y)\in {\mathbb{R}} \,\,is\,\, harmonic$, 
\end{center}
and such that the following conditions are satisfied:
\begin{itemize}
\item[(i)] $w\in {\mathcal{L}}_{\mathrm{loc}}^1(\Omega )$,
\item[(ii)] for each $j$, $1\leq j\leq p$, 
{\Large{$\frac{\partial ^2w}{\partial x_j^2}$}}$\in {\mathcal{L}}_{{\mathrm{loc}}}^1
(\Omega )$ and for each $k$, $1\leq k\leq q$, 
{\Large{$\frac{\partial ^2w}{\partial y_k^2}$}}$\in {\mathcal{L}}_{{\mathrm{loc}}}^1
(\Omega )$.
\end{itemize}
Then $w$ has a separately  harmonic extension to $\Omega$.}}

\vskip0.3cm

\noindent{\textit{Proof.}} Follows directly from Theorem~4 and from the above Lemma of  Federer.
\hfill \qed
\section{Extension results for holomorphic functions}
\subsection{}
Then our result for holomorphic functions, also which is based on Theorem~1. See also our preliminary result \cite{Ri16$_2$}, Theorem~3, p.~156.

\vskip0.3cm

\noindent \textbf{Theorem~5.}
{\emph{Suppose that  $\Omega$ is a domain in ${\mathbb{C}}^n$, $n\geq 1$. 
Let $E\subset \Omega$ be closed in $\Omega$  
and let ${\mathcal{H}}^{2n-1}(E)<+\infty$. Let 
$f=u+iv:\Omega \setminus E\to {\mathbb{C}}$ be holomorphic and such that the following conditions are satisfied:
\begin{itemize}
\item[(i)] $f\in {\mathcal{L}}_{{\mathrm{loc}}}^1(\Omega )$.
\item[(ii)] For each $j$, $1\leq j\leq 2n$, $\frac{\partial ^2u}{\partial x_j^2}\in {\mathcal{L}}_{{\mathrm{loc}}}^1(\Omega )$ and $\frac{\partial ^2v}{\partial x_j^2}\in {\mathcal{L}}_{{\mathrm{loc}}}^1(\Omega )$.
\item[(iii)] For each $j$, $1\leq j\leq 2n$,  and for ${\mathcal{H}}^{2n-1}$-almost all $X_j\in {\mathbb{R}}^{2n-1}$   such that $E(X_j)$ is finite, the following condition holds:\newline
For each $x_j^0\in E(X_j)$ there exist sequences $x_{j,l}^{0,1}, x_{j,l}^{0,2}\in (\Omega \setminus E)(X_j)$, $l=1,2,\dots$, such that $x_{j,l}^{0,1}\nearrow x_{j}^{0}$, $x_{j,l}^{0,2}\searrow x_{j}^{0}$ as $l\to +\infty$, and 
\begin{center}
\begin{itemize}
\item[(iii(a))] $\lim_{l\to +\infty}f(x_{j,l}^{0,1},X_j)=\lim_{l\to +\infty}f(x_{j,l}^{0,2},X_j)\in {\mathbb{C}},$
\item[(iii(b))] $-\infty <\liminf_{l\to +\infty}\frac{\partial u}{\partial x_j}(x_{j,l}^{0,1},X_j)= \limsup_{l\to +\infty} 
\frac{\partial u}{\partial x_j}(x_{j,l}^{0,2},X_j)<+\infty $ \,
and \, $-\infty <\liminf_{l\to +\infty}\frac{\partial v}{\partial x_j}(x_{j,l}^{0,1},X_j)= \limsup_{l\to +\infty}
\frac{\partial v}{\partial x_j}(x_{j,l}^{0,2},X_j)<+\infty $.
\end{itemize}
\end{center}
 \end{itemize}
Then $f$ has a unique holomorphic  extension to $\Omega$.}}
\vskip0.3cm
 \noindent{\textit{Proof.}} It is sufficient to show that $u$ and $v$ have harmonic  extensions $u^*$ and $v^*$ to $\Omega$. As a matter of fact, then $f^*=u^*+iv^*:\Omega \to {\mathbb{C}}$ is   in ${\mathcal{C}}^{\infty}(\Omega)$ and thus a  continuous function. Therefore the claim follows  from Shiffman's theorem or also from  \cite{{Ri78_1},{Ri78_2}}.
 
Another possibility for the proof is just to observe that the in $\Omega \setminus E$  harmonic functions $u$ and $v$ have by Theorem~3 harmonic extensions $u^*$ and $v^*$ to $\Omega$. Since $u^*$ and $v^*$   are thus  ${\mathcal{C}}^{\infty}$ functions, the holomorphy of the  extension $f^*=u^*+iv^*$ in $\Omega$ follows easily. \hfill \qed
\subsection{ } As a concise corollary we get the following: 
\vskip0.3cm
\noindent \textbf{Corollary~5.} (\cite{Ri16$_1$}, Theorem~3, p.~51, \cite{Ri16$_2$}, Theorem~4, p.~157)
{\emph{Suppose that  $\Omega$ is a domain in ${\mathbb{C}}^n$, $n\geq 1$. 
Let $E\subset \Omega$ be closed in $\Omega$ 
and let ${\mathcal{H}}^{2n-1}(E)= 0$. Let 
$f:\Omega\setminus E \to {\mathbb{C}}$ be holomorphic and such that the following conditions are satisfied:
\begin{itemize}
\item[(i)] $f\in {\mathcal{L}}_{{\mathrm{loc}}}^1(\Omega )$,
\item[(ii)] for each $j$, $1\leq j\leq 2n$, $\frac{\partial ^2f}{\partial x_j^2}\in {\mathcal{L}}_{{\mathrm{loc}}}^1(\Omega )$.
 \end{itemize}
Then $f$ has a unique holomorphic extension  to $\Omega$.}}
\subsection{ } Observe that our Theorem~5 and Corollary~5 can be considered, at least in some sense, as  counterparts to  two of Shiffman's well-known extension results for holomorphic functions, namely the following results:

\vskip0.3cm

\noindent \textbf{Shiffman's theorem.} (\cite{Sh68}, Lemma~3, p.~115, and \cite{HaPo75}, Theorem~1.1 (b), p.~703) {\emph{Let
$\Omega$ be a domain in ${\mathbb{C}}^n$, $n\geq 1$. Let $E\subset \Omega$ be closed in $\Omega$ and let ${\mathcal{H}}^{2n-1}(E)<+\infty$. If  $f: \Omega  \rightarrow {\mathbb{C}}$ is  continuous and $f\vert \Omega \setminus E$ is holomorphic, then $f$ is holomorphic in $\Omega$.}}

\vskip0.3cm

\noindent \textbf{Another theorem of Shiffman.} (\cite{Sh68}, Lemma~3, p.~115, and \cite{HaPo75}, Theorem~1.1 (c), p.~703) {\emph{Let
$\Omega$ be a domain in ${\mathbb{C}}^n$, $n\geq 1$. Let $E\subset \Omega$ be closed in $\Omega$ and let ${\mathcal{H}}^{2n-1}(E)=0$. If  $f: \Omega \setminus E \rightarrow {\mathbb{C}}$ is  holomorphic and bounded, then $f$ has a  unique holomorphic extension to $\Omega$.}}

\vskip0.3cm

Shiffman's proofs of his above results were based on coordinate rotation, on the use of Cauchy integral formula, on the already stated important result of Federer, the Lemma above, and on the following classical result of Besicovitch:

\vskip0.3cm

\noindent \textbf{Besicovitch's theorem.} (\cite{Be31}, Theorem~1, p.~2) {\emph{Let
$D$ be a domain in ${\mathbb{C}}$. Let $E\subset D$ be closed in $D$ and let ${\mathcal{H}}^1(E)=0$. If  $f: D\setminus E \rightarrow {\mathbb{C}}$ is  holomorphic and bounded, then $f$ has a unique holomorphic extension to $D$.}}

\vskip0.3cm

For slightly more general versions of Shiffman's results with different proofs, see \cite{Ri78_1}, Theorem 3.1, p.~49, Theorem~3.5, p.~52, Corollary~3.7, p.~54, and \cite{Ri78_2}, Theorem~3.1, p.~333, Corollary~3.3, p.~336. 
\subsection{A previous, slightly related result} Observe that, in addition to Corollary~5, also the following result holds:

\vskip0.3cm

\noindent \textbf{Theorem~6.}
{\emph{Suppose that  $\Omega$ is a domain in ${\mathbb{C}}^n$, $n\geq 1$. 
Let $E\subset \Omega$ be closed in $\Omega$ 
and let ${\mathcal{H}}^{2n-1}(E)= 0$. Let 
$f:\Omega\setminus E \to {\mathbb{C}}$ be holomorphic. If for each $j$, $1\leq j\leq 2n$, $\frac{\partial f}{\partial x_j}\in {\mathcal{L}}_{{\mathrm{loc}}}^2(\Omega )$, then $f$ has a holomorphic extension  to $\Omega$.}}

\vskip0.3cm

The proof follows at once from the following, rather old result:

\vskip0.3cm

\noindent \textbf{Proposition~2.} (\cite{HyRi87}, Corollary~3.6, p. 301)
{\emph{Suppose that  $\Omega$ is a domain in ${\mathbb{C}}^n$, $n\geq 1$. 
Let $E\subset \Omega$ be closed in $\Omega$ 
and let ${\mathcal{H}}^{2n-1}(E)= 0$. Let 
$f:\Omega\setminus E \to {\mathbb{C}}$ be holomorphic. If for some $p\in {\mathbb{R}}$,
\[\int_{\Omega\setminus E}\vert f(z)\vert^{p-2}\sum_{j=1}^n \vert \frac{\partial f}{\partial z_j}(z)\vert^2 \, dm_{2n}(z)<+\infty ,\]
then $f$ has a meromorphic extension $f^*$ to $\Omega$. If $p\geq 0$, then $f^*$ is holomorphic.}}

\vskip0.3cm

For related, partly previous and partly more general results, see \cite{Ce83}, Theorem, p.~284, \cite{HyRi87}, Theorem~3.5, pp.~300-301, and \cite{Ri90}, Theorem~3.1, pp.~925-926.


\vskip1cm

\begin{thebibliography}{99}
\footnotesize{
\bibitem{ArGa01}  Armitage, David H., Gardiner, Stephen J.: Classical Potential Theory.  Springer Monographs in Mathematics, Springer-Verlag, London (2001)
\bibitem{Be31} Besicovitch, A.S.: On sufficient conditions for a function to be analytic, and on behavior of analytic functions in the neighborhood of non-isolated singular point. Proc. London Math. Soc. (2) {\textbf{32}}, 1-9 (1931)
\bibitem{Bl95} Blanchet, P.: On removable singularities of subharmonic and plurisubharmonic functions.
Complex Variables {\textbf{26}}, 311-322 (1995)
\bibitem{Ce83} Cegrell, U.: Removable singularity sets for analytic functions having modulus with bounded Laplace mass. Proc. Amer. Math. Soc. {\textbf{88}}, 283-286 (1983)
\bibitem{Ch89} Chirka, E.M.: Complex Analytic Sets. Kluwer Academic Publisher, Dordrecht (1989)
\bibitem{EvGa93} Evans, Lawrence C., Gariepy, Ronald F.: Measure Theory and Fine Properties of Functions. Studies in Advanced Mathematics, CRC Press, Inc,, Boca Raton (1992)
\bibitem{Fe69} Federer, H.: Geometric measure theory. Springer, Berlin (1969)
\bibitem{HaPo70} Harvey, R., Polking J.: Removable singularities of solutions of linear partial differential equations. Acta Math. {\textbf{125}}, 39-56 (1970)
\bibitem{HaPo75} Harvey, R., Polking, J.: Extending analytic objects. Comm. Pure and Appl. Mathematics {\textbf{28}}, 701-727 (1975)
\bibitem{HaKe76} Hayman, W.K., Kennedy, P.B.: Subharmonic Functions, Vol. I.  Academic Press, London (1976)
\bibitem{Hel69} Helms,  L.L.: Introduction to potential theory.  Wiley-Interscience, New York (1969)
\bibitem{Her71}  Hervé, M.: Analytic and plurisubharmonic functions in finite and infinite
dimensional spaces.  Lecture Notes in Mathematics, vol. 198. Springer, Berlin (1971)
\bibitem{HyRi87} Hyvönen, Jaakko, Riihentaus, Juhani: Removable singularities for holomorphic functions with locally finite Riesz mass. J. London Math. Soc. (2) {\textbf{35}}, 296-302 (1987)
\bibitem{JaPf00} Jarnicki, M., Pflug, P.: Extension of Holomorphic Functions.  Walter de Gruyter, Berlin (2000)
\bibitem{JaPf11} Jarnicki, M., Pflug, P.: Separately Analytic Functions.  European Mathematical Society, Z\"urich (2011)
\bibitem{Kh92}  Khabibullin, B.N.: A uniqueness theorem for subharmonic functions of finite order. Mat. Sb. {\textbf{182}}(6), 811-827 (1991); English transl. in Math USSR Sbornik {\textbf{73}}(1), 195-210 (1992)
\bibitem{Kh99}  Khabibullin, B.N.: Completeness of systems of entire functions in spaces of holomorphic functions. Mat. Zametki {\textbf{66}}(4), 603-616 (1999); English transl. in Math. Notes {\textbf{66}}(4), 495-506 (1999)
\bibitem{Le61} Lelong, P.: Fonctions plurisousharmoniques et fonctions analytiques réelles. Ann. Inst. Fourier {\textbf{11}}, 515-562 (1961)
\bibitem{Le69} Lelong, P.: Plurisubharmonic functions and positive differential forms.  Gordon and Breach, New York (1969)
\bibitem{Po84} Polking, J.: A survey of removable singularities. Seminar on nonlinear partial differential equations, Berkeley, California, 1983. In: Math. Sci. Res. Inst. Publ (ed. S.S. Chern) {\textbf{2}}, 261-292 (1984), Springer, Berlin (1984)
\bibitem{Ra37} Rado,  T.: Subharmonic functions.  Springer, Berlin (1937)
\bibitem{Ri78_1}  Riihentaus, J.: Removable singularities of analytic functions of several complex variables. Math. Z. {\textbf{32}}, 45-54 (1978)
\bibitem{Ri78_2} Riihentaus, J.: Removable singularities of analytic and meromorphic functions of
several complex variables.  Colloquium on Complex Analysis, Joensuu, Finland, August 24-27, 1978 
(Complex Analysis, Joensuu 1978). In: Proceedings (eds. Ilpo Laine, Olli Lehto, Tuomas Sorvali),
Lecture Notes in Mathematics, vol. 747, 329-342 (1978), Springer, Berlin (1979)
\bibitem{Ri90} Riihentaus, Juhani: A nullset for normal functions in several variables. Proc. Amer. Math. Soc. {\textbf{110}}(4), 923-933 (1990)
\bibitem{Ri04} Riihentaus, J.: Subharmonic functions, mean value inequality, boundary behavior, nonintegrability and exceptional sets.
Workshop on Potential Theory and Free Boundary Flows; August 19-27, 2003: Kiev, Ukraine. In: Transactions of the Institute of Mathematics of the National Academy of Sciences of Ukraine, Kiev {\textbf{1}}(3), 169-191 (2004)
\bibitem{Ri07} Riihentaus, J.: Subharmonic functions, generalizations and separately subharmonic functions. The XIV-th Conference on Analytic Functions,
Che\l m, Poland, July 22-28, 2007 (eds. Dariusz Partyka and Józef Zajac). In: Scientific Bulletin of Che\l m, Section of Mathematics and Computer Science {\textbf{2}}, 49-76 (2007) (ISBN 978-83-61149-24-8) (arXiv:math/0610259v5 8 Oct 2008)}
\bibitem{Ri09} Riihentaus, J.: Subharmonic functions, generalizationas, weighted boundary behavior, and separately subharmonic functions: A survey. Fifth World Congress of Nonlinear Analysts (WCNA 2008), Orlando, Florida, USA,  July 2-9, 2008. In: Nonlinear Analysis, Series A: Theory, Methods \& Applications {\textbf{71}}(12), e2613-e2627 (2009) (doi: 10.1016/j.na.2009.05.077) 
\bibitem{Ri13} Riihentaus, J.: An inequality type condition for quasinearly subharmonic functions and applications. Positivity VII, Leiden, July 22-26, 2013, Zaanen Centennial Conference. In: Ordered Structures and Applications: Positivity VII, Trends in Mathematics, 395-414, Springer International Publishing (2016)
\bibitem{Ri15} Riihentaus, J.: Exceptional sets for subharmonic functions. J. Basic \& Applied Sciences  {\textbf{11}}, 567-571 (2015)
\bibitem{Ri16$_1$} Riihentaus, J.: A removability result for holomorphic functions of several complex variables. J. Basic \& Applied Sciences  {\textbf{12}}, 50-52 (2016)
\bibitem{Ri16$_2$} Riihentaus, J.: Removability results for subharmonic functions, for harmonic functions and for holomorphic functions. Matematychni Studii {\textbf{46}}(2), 152-158 (2016)
\bibitem{Ri19} Riihentaus, J.: A removability result for separately subharmonic functions. Visnyk of the Lviv Univ. Series Mech. Math. 2018. Issue 85, 60-65 (2019)
\bibitem{Sh68} Shiffman, B.: On the removal of singularities of analytic sets. Michigan Math. J. {\textbf{15}}, 111-120 (1968)
\bibitem{Ul91} Ullrich, D.C.: Removable sets for harmonic functions. Michigan Math. J. {\textbf{38}}, 467-473 (1991)

\end{thebibliography}
\end{document}